\documentclass[reqno]{amsart}
\usepackage{hyperref}
\usepackage{mathrsfs}
\usepackage{mathtools}

\begin{document}
\title[]
{weak solutions to the time-fractional $g-$B\'enard equations}

\author[K. Akhlil]
{Khalid Akhlil} 

\author[S. Ben Aadi]
{Sultana Ben Aadi}

\author[H. Mahdioui]
{Hicham Mahdioui}

\email{k.akhlil@uiz.ac.ma}
\email{sultana.benaadi@edu.uiz.ac.ma}
\email{h.mahdioui@uiz.ac.ma}

\date{\today}
\thanks{}

\begin{abstract}
In this paper, we introduce the $g$-B\'enard equations with time-fractional derivative of order $\alpha\in(0,1)$ in domains of $\mathbb R^2$. This equations model, the memory-dependent heat conduction of liquids in fractal media considered in $g$-framework. We aim to study the existence and uniqueness of weak solutions by means of standard techniques from Navier-Stokes equations theory and fractional calculus theory.

\end{abstract}

\maketitle \setlength{\textheight}{19.5 cm}
\setlength{\textwidth}{12.5 cm}
\newtheorem{theorem}{Theorem}[section]
\newtheorem{lemma}[theorem]{Lemma}
\newtheorem{proposition}[theorem]{Proposition}
\newtheorem{corollary}[theorem]{Corollary}
\newtheorem{Hypo}[theorem]{Assumption}
\theoremstyle{definition}
\newtheorem{definition}[theorem]{Definition}
\newtheorem{assumptions}[theorem]{Assumptions}
\newtheorem{example}[theorem]{Example}
\theoremstyle{remark}
\newtheorem{remark}[theorem]{Remark}
\numberwithin{equation}{section} \setcounter{page}{1}

\section{ Introduction }
Let $\Omega_g=\Omega_2\times (0,g)$ where $\Omega_2$ is a bounded domain in $\mathbb R^2$ and $g$ is some scalar nonnegative function. We introduce the following time-fractional $g$-B\'enard equations of the following form:
  \begin{align*}
  \partial^{\alpha}_t u+(u\cdot\nabla)\,u-\nu  \Delta u+\nabla p&= \xi \theta +f_1(t)  \\[0.1cm]
  \nabla\cdot gu &= 0\\[0.1cm]
  \partial^{\alpha}_t \theta+(u\cdot\nabla)\,\theta-\kappa\Delta\theta &= f_2(t)
  \end{align*}
where $u$ is the fluid velocity, $p$ is the pressure, $\theta$ is the temperature, $f_1$ is the external force function, $f_2$ is the heat source function, $\xi\in\mathbb{ R}^3$ is a constant vector, $\nu$ the kinematic viscosity and $\kappa$ the thermal diffusivity are positive constants. The derivative of order $\alpha$ is considered in the Caputo sense. The time-fractional $g$-B\'enard problem consists in a system that couples time-fractional Navier-Stokes equations and time-fractional advection-diffusion heat equation in order to model a memory-dependent convection in a fluid considered in a fractal media.

The introduction of time-fractional derivative in fluid dynamics goes back to Lions in \cite{L59} but for order less than $\frac{1}{4}$ provided the space dimension is not further than $4$.  In recent works of Zhou and Peng \cite{ZP17, ZP172}, the question of weak solutions and optimal control problem of time-fractional Navier-Stokes equations in fractal media was considered. Numerical results regarding such problems was treated firstly in \cite{KKAR14} and constitute an emerging field of research. More recently, time-fractional $g$-Navier-Stokes problem was introduced and results regarding existence, uniqueness of solutions and optimal control are proved\cite{BA20}. 

The theory of $g$-Navier-Stokes equations started with the works of Hale and Raugel \cite{HR92a, HR92b}, Raugel and Sell \cite{RS93} who studied 3d nonlinear equations and Navier-Stokes equations in thin domains. J. Roh\cite{R01}, a student of Sell, generalized the previous works to thin domains of the form $ \Omega_g=\Omega\times(0,g)$, where $g$ is some smooth scalar function. The derived equations are called the $g-$Navier-Stokes equations. This theory have interested many researchers in recent years, see  \cite{BR04,AQ12a,JHW11} and references therein. 

On the other hand, heat conduction based on the classical Fourier law, which relates the heat flux vector and the temperature gradient, has shown its limits. The time-fractional heat conduction model can be seen as a good alternative\cite{MLP01,GIL00,SZ97}, and references therein. Boussinesq(or B\'enard) model is a combination of the heat conduction model and Navier-Stokes equations and is a well developed subject in modelling heat conducting fluids\cite{Z89,Z03,GLCS80,CD80}. The aim of this paper is to generalize the setting in \cite{OK16,OK17}, where $g$-B\'enard equations was considered, to time-fractional $g$-B\'enard equations.

This paper is organized as follows: In the first section, we will recall some concept and notations related to fractional calculus. Section 2, is devoted to the problem statement. Section 3 will be dedicated to the proof of the existence and uniqueness of weak solutions to time-fractional $g$-Navier-Stokes equations.

\section{Preliminaries on fractional calculus}

In this section, we provide some notations and preliminary results concerning fractional calculus. For this purpose, assume $X$ to be a Banach space. Let $\alpha\in (0,1]$ and let $k_\alpha$ denote the Riemann-Liouville kernel
\[
k_\alpha(t)=\frac{t^{\alpha-1}}{\Gamma(\alpha)}
\]

For a function $v:[0,T]\rightarrow X$, we give the following definitions of derivatives and integrals:
\begin{itemize}
\item[(1)] The left Riemann-Liouville integral of $v$ is defined by
\[
I_t^\alpha v(t)=\int_0^t k_\alpha(t-s)v(s)\,ds,\quad t>0
\]provided the integral is point-wise defined on $[0,+\infty[$
\item[(2)] The right Riemann-Liouville integral of $v$ is defined by
\[
\, I_{t,T}^\alpha v(t)=\int_t^T k_\alpha(t-s)v(s)\,ds,\quad t>0
\]provided the integral is point-wise defined on $[0,+\infty[$
\item[(3)]  The left Caputo fractional derivative of order $\alpha$ of $v$, is defined by
\[
\,D_t^\alpha v(t)= \int_0^t k_{1-\alpha}(t-s)\frac{d}{ds}v(s)\,ds
\]
\item[(4)] The right Riemann-Liouville fractional derivative of order $\alpha$ of $v$ is defined by
\[
\,D_{t,T}^\alpha v(t)=-\frac{d}{dt} \int_t^T k_{1-\alpha}(t-s)v(s)\,ds
\]
\item[(5)] The Liouville-Weyl fractional integral on the real axis for functions $v:\mathbb R\rightarrow X$ is defined as follows
\[
\,I_{-,t}^\alpha v(t)=\int_{-\infty}^t k_{\alpha}(t-s)v(s)\,ds
\]
\item[(6)] The Caputo fractional derivative on the real axis for functions $v:\mathbb R\rightarrow X$ is defined as follows
\[
\,D_{-,t}^\alpha v(t)= \,I_{-,t}^{1-\alpha}\frac{d}{dt} v(t)
\]
\end{itemize}

Note that the notation $\partial_t^\alpha$ stands for Caputo fractional partial derivative, i.e. when functions have another argument than time. We have the following fractional integration by parts formula; see, e.g. \cite{A07}
\begin{align}\label{partial}
\int_0^T(\partial_t^\alpha u(t),\psi(t))dt&=\int_0^T(u(t),D_{t,T}^\alpha\psi(t))dt+(u(t),I_{t,T}^{1-\alpha}\psi(t))|_0^T\\\nonumber
&=\int_0^T(u(t),D_{t,T}^\alpha\psi(t))dt-(u(0),I_T^{1-\alpha}\psi(t))
\end{align}since for $\psi\in C_0^\infty([0,T],X)$ one have $\displaystyle\lim_{t\to T}\,I_{t,T}^{1-\alpha}\psi(t)=0$.

To pass from weak convergence to strong convergence we will need a compactness result. Let $X_0$, $X$, $X_1$ be Hilbert spaces with $ X_0 \xhookrightarrow{} X \xhookrightarrow{} X_1$ being continuous and $X_0\xhookrightarrow{} X$ being compact. Assume that $v:\mathbb R\rightarrow X_1$ and denote by $\widehat v$ its Fourier transform
\[
\widehat v(\tau)=\int_{-\infty}^{+\infty}e^{-2i\pi t \tau}v(t)\,dt
\]
We have for $\gamma>0$
\[
\widehat{D_t^\gamma v}(\tau)=(2i\pi\tau)^\gamma \widehat v(\tau)
\]
For a given $0<\gamma<1$, we introduce the following space
\[
W^\gamma(\mathbb R,X_0,X_1)=\left\{v\in L^2(\mathbb R,X_0):\, D_t^\gamma v\in L^2(\mathbb R,X_1)\right\}
\]
Clearly, it is a Hilbert space for the norm
\[
\|v\|_\gamma=\left(\|v\|^2_{L^2(\mathbb R,X_0)}+\|\,|\tau|^\gamma\widehat v\|^2_{L^2(\mathbb R,X_1)}\right)^{1/2}
\]
 For any set $K\subset \mathbb R$, we associate with it the subspace $W_K^\gamma\subset W^\gamma$ defined as
 \[
 W_K^\gamma(\mathbb R,X_0,X_1)=\{v\in W^\gamma(\mathbb R,X_0,X_1):\,\mathrm{support }\, u\subset K\}
 \]

By similar discussion as in the proof of Theorem 2.2 in Temam \cite{T84}(see also Theorem 2.1 in \cite{ZP17}), it is clear that $W_K^\gamma(\mathbb R,X_0,X_1)\xhookrightarrow{} L^2(\mathbb R,X)$ is compact for any bounded set $K$ and any $\gamma>0$.

As a particular situation of the compactness result discussed above, let $H$, $V$ be two Hilbert spaces endowed with the scalar product $(.,.)_H$ and $(.,.)_V$ and the norms $|.|_H$ and $\|.\|_V$, respectively. Denote by $\langle.,.\rangle$ the dual pairing between $V$ and $V'$, the dual of $V$. Moreover assume that $V \xhookrightarrow{} H  \xhookrightarrow{} V'$ continuously and compactly and note that the space
\[
W^\gamma(0,T;V,V')=\left\{v\in L^2(0,T;V):\, \partial_t^\gamma v\in L^2(0,T;V')\right\}
\]is compactly embedded in $L^2(0,T;H)$. Similarly to Lemma 2.1 in \cite{ZP17}, we have
\[
\partial_t^\gamma(u(t),v)_V=\langle\partial_t^\gamma u(t),v\rangle
\] for $u\in W^\gamma(0,T;V,V')$ and $v\in H$. Moreover, for a derivable function $v:[0,T]\rightarrow V$ we have from \cite{A10} that

\[
(v(t),D_t^\gamma v(t))_{H}\geq\frac{1}{2}\,D_t^\gamma|v(t)|^2
\]

We end this section by the following important result
\begin{lemma}

Suppose that a nonnegative function satisfies
\[
_{~0}^{~C}D_t^\gamma v(t)+ c_1 v(t)\leq c_2(t)
\]for $c_1>0$ and $c_2$ a nonnegative integrable function for $t\in[0,T]$. Then 
\[
v(t)\leq v(0)+\frac{1}{\Gamma(\gamma)}\int_0^t(t-s)^{\gamma-1} c_2(s)\,ds
\]

\end{lemma}

For more details about fractional calculus we refer to the monographs \cite{KST06,Z14,Z16}.

\section{ Problem Statement}

We introduce the usual notation used in the context of the mathematical theory of Navier-Stokes equations \cite{T84}. Let $\Omega_g = \Omega_2 \times (0, g) = (0, 1)\times (0, 1)\times (0, g)$ where
$g = g(y_1, y_2)$ is a smooth function defined on $\Omega_2$. In addition we assume that

\begin{equation}\label{2.1}
  \begin{array}{c }
\displaystyle  0<m_0< g(y_1,y_2)\leq M_0,~~~~~~\text{for all}~~ (y_1, y_2)\in \Omega_2
  \\ [0.2cm]
  \displaystyle |\nabla g|_{\infty}=\displaystyle\sup_{\Omega_2}| \nabla g  |<\infty,~~~~~~~~g\in C^{\infty}_{per}(\Omega_2).
\end{array}
\end{equation}

Let $L^2(\Omega, g)$ denotes the Hilbert space, of weighted Sobolev spaces type, with the inner product \[\langle u,v\rangle_g = \int_\Omega(u \cdot v)gdx\] and the induced norm $|u|^2_g = \langle u, u\rangle_g$. Similarly, we can define the weithed Sobolev space $H^1(\Omega, g)$ equipped with the norm

\[ |u|^2_{H^1(\Omega,g)}=\langle u,u \rangle_g + \displaystyle \sum^n_{i=1} \langle \frac{\partial u}{\partial x_i} ,\frac{\partial u}{\partial x_i} \rangle_g.  \]
Moreover, we will need the following spaces:

\begin{align*}
\mathcal{V}_1 & = \, \{ u\in (C^{\infty}_{per}(\Omega))^n~~: ~~ \nabla\cdot(gu)=0,~~~\displaystyle\int_{\Omega}udx=0~~~\text{on}~~~\Omega      \}\\
H_g & =  \text{ the closure of } \mathcal{V}_1 \text{ on}~~~ L^2(\Omega,g)\\
 V_g & =   \text{ the closure of }\mathcal{V}_1 \text{ on }~~~ H^1(\Omega,g) \\
 V'_g & =    \text{  the dual space of  }~~~V_g \\
  \mathcal{V}_2 & = \, \{ \varphi\in C^{\infty}_{per}(\Omega)~~:~~ \displaystyle\int_{\Omega}\varphi dx=0    \}\\
 W_g & =  \text{ the closure of } \mathcal{V}_2 \text{ on}~~~ H^1(\Omega,g) \\ 
  W'_g & =  \text{  the dual space of  }~~~W_g\\
    Q & =  \text{  the closure of}~~~ \{\nabla\varphi~~:~~\varphi\in C^1_{per}(\overline{\Omega},R)  \} ~~\text{in}~~L^2(\Omega)
\end{align*}
where $H_g$ is endowed with the inner product and the norm in $L^2(\Omega g)$. In addition, the spaces  $V_g$ and $Wg$ are endowed with the inner product and the norm in $H^1(\Omega, g)$. Let us also remark that the inclusions
\[V_g\subset  H_g = H'_g\subset V'_g,\] \[W_g \subset L^2(\Omega, g) \subset W'_g,\]
are dense and continuous \cite{M92,R01}. By the Riesz representation theorem, it is possible to write
\[ \langle f,u \rangle_g=(f,u)_g,~~~~~~~~ \forall f\in H_g,~\forall u\in V_g.      \]

Let us now define the orthogonal projection $P_g$ as $P_g : L^2_{ per}(\Omega, g)\rightarrow  H_g$. It is clear that $Q \subseteq H^{\perp}_{g}$ . Similarly, we define $\widetilde{P}_g$ as $\widetilde{P}_g :L^2_{ per}(\Omega, g)\rightarrow  W_g$. By taking into account the following equality\cite{R01}:
\[  \displaystyle -\frac{1}{g} ( \nabla\cdot g\nabla u )  =-\Delta u - \frac{1}{g}( \nabla g\cdot \nabla  )u\]we define the $g-$Laplace operator and $g$-Stokes operator as follows
\[-\Delta_gu=-\displaystyle \frac{1}{g}( \nabla\cdot g\nabla u)  \] and
\[A_gu = P_g[-\Delta_gu]\] respectively. We have the following result \cite{R01}:

\begin{proposition} For the g-Stokes operator $A_g$, the followings hold:

\begin{enumerate}
  \item[(1)] The g-Stokes operator $A_g$ is a positive, self-adjoint operator with compact inverse, where the domain of $A_g$ is
  $D(A_g) = V_g \cap H^2(\Omega, g).$
  \item[(2)]  There exist countable eigenvalues of $A_g$ satisfying
  \[     0< \displaystyle \frac{4\pi^2m_0}{M_0}\leq  \lambda_1  \leq \lambda_2 \leq\lambda_3\leq\cdot\cdot\cdot   \]
where $\lambda_1$ is the smallest eigenvalue of $A_g$. In addition, there exist the corresponding collection of eigenfunctions $\{w_i\}_{i\in\mathbb{N}} $ forms an orthonormal basis for
$H_g$.
\end{enumerate}
\end{proposition}

The operators $A_g$ and $P_g$ are clearly self-adjoint, then by using integration by parts we have:
  \begin{align*} 
  \langle A_g u,u\rangle_g= &\langle P_g[ \displaystyle -\frac{1}{g}(\nabla\cdot g \nabla )u],u\rangle_g\\
    =&\displaystyle \int_{\Omega}  (\nabla u\cdot\nabla u )_gdx \\
    = &  \langle \nabla u\cdot\nabla u \rangle_g
 \end{align*}
It then follows that for $u \in V_g$ we can write $  | A^{1/2}u|_g=  | \nabla u |_g=\| u \|_g   $. On the other hand,, since the functional
\[  \tau\in W_g\rightarrow ( \nabla \theta,\nabla\tau)_g\in \mathbb{R}    \]
is a continuous linear mapping on $W_g$, we can define a continuous linear mapping $\widetilde{A}_g$ on $W'_g$ such that
\[    \forall \tau \in W_g,~~~~~\langle \widetilde{A}_g,\tau \rangle_g = ( \nabla \theta  , \nabla \tau )_g \]
for all $\theta\in  W_g$. For $u$, $v$ and $w$ laying in an appropriate subspaces of $L^2_{per}(\Omega, g)$, we can define the bilinear operator \[B_g(u, v) = P_g[(u \cdot \nabla)v]\] and the trilinear form
\[  b_g(u,v,w)=\displaystyle\sum^{n}_{i,j=1}\int_{\Omega}   u_i\frac{\partial v_j}{\partial x_i} w_j gdx \]
As a consequence, one obtains $b_g(u, v, w) =
-b_g(u, w, v)$ which implies that $b_g(u, v,v) = 0$. MOreover, we have the following inequality satisfied by $b_g$(see \cite{T84,W09}),
\begin{equation}\label{2.7}
| b_g(u,v,w)|_g\leq   c| u|^{1/2}_{g}  \| u\|^{1/2}_{g}  | v|_{g}   | w|^{1/2}_{g}  \| w\|^{1/2}_{g},~~~~~~\forall u,v,w\in V_g
\end{equation}

Similarly, for $u \in V_g$ and $\theta, \tau\in W_g$ we define $\widetilde{B}_g(u,\theta) = \widetilde{P}_g[(u \cdot \nabla)\theta]$ and
\[  \widetilde{b}_g(u,\theta,\tau)=\displaystyle\sum^{n}_{i,j=1}\int_{\Omega}   u_i(x)\frac{\partial \theta(x)}{\partial x_j} \tau(x) gdx \]

We denote the operators
$  C_gu = P_g[\frac{1}{g} (\nabla g \cdot \nabla)u]$  and $\widetilde{C}_g \theta = \widetilde{P}_g[\frac{1}{g} (\nabla g \cdot \nabla)\theta]$ such that
\begin{align*}
 \langle  C_gu,v \rangle_g  =&b_g (\displaystyle \frac{\nabla g}{g},u,v),\\
 \langle \widetilde{C}_g \theta,\tau  \rangle_g=&\widetilde{b}_g (\displaystyle \frac{\nabla g}{g},\theta,\tau  )
\end{align*}
Finally, let $\widetilde{D}_g\theta= \widetilde{P}_g[\displaystyle \frac{\nabla g}{g}\theta ]$ such that
\[ \langle  \widetilde{D}_g\theta ,\tau  \rangle_g=-\widetilde{b}_g (\displaystyle \frac{\nabla g}{g},\theta,\tau  ) - \widetilde{b}_g (\displaystyle \frac{\nabla g}{g},\theta,\tau  ). \]
We can now rewrite the system of $g$-B\'enard equations in the following abstract time-fractional evolutionary equations
\begin{equation}\label{benard}
\begin{array}{ccc}
&\partial^{\alpha}_{t}u+ B_g(u,u)+\nu A_gu+\nu C_gu= \xi\theta+f_1\\[0.2cm]
&\partial^{\alpha}_{t} \theta+ \widetilde{B}_g(u,\theta)+\kappa \widetilde{A}_g\theta-\kappa \widetilde{C}_g\theta-   \kappa\widetilde{D}_g\theta= f_2\\[0.2cm]
&u(x,0)=u_0(x),~~~~~~~~\theta(x,0)=\theta_0(x)
\end{array}
\end{equation}

We give the following two lemmas and for the proofs of these lemmas we refer to \cite{BR04,L69}.

\begin{lemma}\label{lem1} For $n = 2$ there exists a positive constant $c$ such that

\[ |u|_{L^4(\Omega,g)} \leq c|u|^{1/2}_{g}  |\nabla u|^{1/2}_{g} ,~~~~~~~~ \forall u\in H^1( \Omega , g ) . \]
\end{lemma}

\begin{lemma}\label{lem2} For  $   u\in L^2(  0,T,V_g)$ we have:
\[B_g(u,u)(t)\in L^1(0,T,V'_g)  \text{  and  } C_gu(t)\in L^2(0,T,H_g).      \]
\end{lemma}

%%%%%%%%%%%%%%%%%%%%%%%%%%%%%%%%%%%%%%%%%%%
%%%%%%%%%%%%%%%%%%%%%%%%%%%%%%%%%%%%%%%%%%%
%

%%%%%%%%%%%%%%%%%%%%%%%%%%%%%%%%%%%%%%%%%%
%%%%%%%%%%%%%%%%%%%%%%%%%%%%%%%%%%%%%%%%%%%

%%%%%%%%%%%%%%%%%%%%%%%%%%%%%%%%%%%%%%%%%%%

\section{    Existence of weak solutions   }

\begin{definition} A pair of functions $ \{ u,\theta\}$ is called a weak solution of the system \eqref{benard} if $u\in  L^2(0, T; V_g)$ and $\theta \in L^2(0, T; W_g)$ satisfy the following equations

\begin{equation}\label{3.1}
\begin{array}{ccc}
&\partial^{\alpha}_{t}(u,v)_g+ b_g(u,u,v)+\nu (\nabla u,\nabla v  )_g+\nu (C_gu,v)_g= (\xi\theta,v)_g+(f_1,v)_g\\[0.2cm]
&\partial^{\alpha}_{t} (\theta,\tau)_g+ \widetilde{b}_g(u,\theta,\tau)+\kappa(  \nabla\theta,\nabla\tau )_g+\kappa \widetilde{b}_g
( \displaystyle \frac{\nabla g}{g},\tau,\theta)=( f_2,v)_g
\end{array}
\end{equation}

for all $v_2 \in V_g$ and $\tau\in W_g.$
\end{definition}
\begin{theorem}If $f_1\in  L^{\frac{2}{\alpha_1}}(0, T; L^2(\Omega, g))$ and $f_2\in  L^{\frac{2}{\alpha_2}}(0, T; L^2(\Omega, g))$ $(\alpha_1,\,\alpha_2<\alpha)$, $u_0 \in H_g$, $\theta_0 \in L^2(\Omega, g)$ and $g$ is a smooth function satisfying the conditions given in (\ref{2.1}) defined on $\Omega_2$ then, there exist a unique weak solution $\{ u,\theta \}  $ of the system \eqref{benard} satisfying the periodic boundary conditions.
\end{theorem}

For the proof of the theorem, we shall use the standard Feado- Galerkin method.

\begin{proof}

Since $V_g$ is separable and $\mathcal{V}_1$ is dense in $V_g$, there exists a sequence $\{u_i\}_{i\in \mathbb{N}}$ which forms a
complete orthonormal system in $H_g$ and a base for $V_g$. Similarly, there exists a sequence
$\{\theta_i\}_{i\in \mathbb{N}}$ which forms a complete orthonormal system in $L^2(\Omega, g)$ and a base for $W_g$. Let $m$ be an arbitrary but fixed positive integer. For each $m$, we define an approximate
solution $\{u^{(m)}(t), \theta^{(m)}(t)\}$ of \eqref{benard} for $1 \leq k \leq m$ and $t \in [0 ,T]$ in the form,

\begin{equation}\label{3.3}
\begin{array}{lr}
 \displaystyle u^{(m)}(t)=\sum^m_{j=1} f_j^{(m)}(t)u_j,~~~~~~~& \displaystyle\theta^{(m)}(t)=\sum^m_{j=1} g_j^{(m)}(t)\theta_j
\end{array}
\end{equation}
and we consider the following approximate problem

\begin{equation}\label{3.4}
\begin{array}{l}
\partial^{\alpha}_{t}(u^{(m)},u_k)_g+ b_g(u^{(m)},u^{(m)},u_k)+\nu ((u^{(m)},u_k ) )_g
+\nu b_g( \displaystyle \frac{\nabla g}{g}   ,u^{(m)},u_k)\\
\qquad= (\xi\theta^{(m)},u_k    )_g+(f_1,u_k)_g,
\end{array}
\end{equation}

\begin{equation}\label{3.5}
\begin{array}{l}
\partial^{\alpha}_{t} (\theta^{(m)} ,\theta_k)_g+ \widetilde{b}_g(u^{(m)},\theta^{(m)},\theta_k)+\kappa(( \theta^{(m)},\theta_k ))_g
+\kappa \widetilde{b}_g( \displaystyle \frac{\nabla g}{g},\theta_k,\theta^{(m)})\\
\qquad=( f_2,\theta_k)_g, 
\end{array}
\end{equation}and
\begin{equation}\label{3.6}
\begin{array}{lr}
 \displaystyle u^{(m)}(0)=u_{m_0}=\sum^m_{j=1} (a_0,u_j)u_j,~~~~~~~& \displaystyle\theta^{(m)}(0)=\theta_{m_0}=\sum^m_{j=1} (\tau_0,\theta_j)\theta_j
\end{array}
\end{equation}
This system forms a nonlinear fractional order system of ordinary differential equations for the functions $f_j^{(m)}(t)$ and $g_j^{(m)}(t)$ and has a maximal solution on some interval $[0, T ]$ ( cf. \cite{BA20}). We multiply (\ref{3.4}) and (\ref{3.5}) by $f_j^{(m)}(t)$ and $g_j^{(m)}(t)$ respectively, and add these equations for $k = 1,\cdot\cdot\cdot, m.$ Taking into account $b_g(u^{(m)}, u^{(m)}, u^{(m)}) = 0$ and $\widetilde{b}g(u^{(m)}, \theta^{(m)},\theta^{(m)}) =0$ we get;

\begin{equation}\label{3.7}
\begin{array}{r}
(D^{\alpha}_{t}u^{(m)},u^{(m)})_g+ \nu\| u^{(m)} (t) \|^2_g +\nu b_g( \displaystyle \frac{\nabla g}{g}   ,u^{(m)}(t),u^{(m)}(t))
\\ = (\xi\theta^{(m)},u^{(m)}(t)    )_g+(f_1,u^{(m)}(t))
\end{array}
\end{equation}
and
\begin{equation}\label{3.8}
\begin{array}{r}
 (D^{\alpha}_{t}\theta^{(m)}(t) ,\theta^{(m)}(t))_g+  \kappa\| \theta^{(m)}(t)  \|^2_g   +\kappa \widetilde{b}_g( \displaystyle \frac{\nabla g}{g},\theta^{(m)}(t),\theta^{(m)}(t))=( f_2,\theta^{(m)}(t))_g
\end{array}
\end{equation}

Using Schwarz and Young inequalities in (\ref{3.7}) and (\ref{3.8})

\begin{align*}
&D^{\alpha}_t  |  u^{(m)}(t) |^2_{g} +\nu \| u^{(m)}(t) \|^2_{g}  \leq   \displaystyle \frac{ M_0 |\xi|^2_{\infty}}{\pi^2 m_0\nu} |  \theta^{(m)}(t)   |^2_g                +\frac{4}{\nu}\| f_1(t) \|^2_{V'_g}+\frac{   2\nu |\nabla g|^2_{\infty}}{m^2_0}|   u^{(m)}(t) |^2_g\\
&D^{\alpha}_t  |  \theta^{(m)}(t) |^2_{g} +\kappa \| \theta^{(m)}(t) \|^2_{g}  \leq    \displaystyle \frac{ 2}{\kappa}\| f_2(t) \|^2_{W'_g}
+\frac{2\kappa|\nabla g|^2_{\infty}}{m^2_0}|\theta^{(m)}(t)|^2_g 
\end{align*}
By using the fact that  $\displaystyle |\nabla g |^2_{\infty}< \frac{\pi^2 m^3_0}{M_0}$ and noting $\nu'=\nu (1- \frac{ M_0|\nabla g|^2_{\infty}}{2\pi^2 m^3_0} )$, $\kappa'=\kappa(1- \frac{ M_0|\nabla g|^2_{\infty}}{2\pi^2 m^3_0}      )$ and $c'= \frac{ M^2_0\|\xi\|^2_{\infty}}{4\pi^4 m^2_0}$, we get the inequalities
\begin{equation}\label{3.9}
\begin{array}{l}
D^{\alpha}_t  |  u^{(m)}(t) |^2_{g} +\nu' \| u^{(m)}(t) \|^2_{g}  \leq   \displaystyle \frac{ c'}{\nu} \|  \theta^{m}(t)   \|^2_g                +\frac{4}{\nu}\| f_1(t) \|^2_{V'_g}
\end{array}
\end{equation}
and
\begin{equation}\label{3.10}
\begin{array}{l}
 D^{\alpha}_t  |  \theta^{(m)}(t) |^2_{g} +\kappa' \| \theta^{(m)}(t) \|^2_{g}  \leq    \displaystyle \frac{ 2}{\kappa}\| f_2(t) \|^2_{W'_g}  \end{array}
\end{equation}
Integrating  \eqref{3.10} from $0$ to $T$, in the fractional sense, we obtain

\begin{align*}
|\theta^{(m)}(t)|_g^2 +&\frac{\kappa'}{\Gamma(\alpha)}\int_0^t (t-s)^{\alpha-1}\|\theta^{(m)}(s)\|_g^2\, ds\\\leq&\, |\theta_{0m}|_g^2+\frac{2}{\kappa\Gamma(\alpha)}\int_0^t(t-s)^{\alpha-1}\| f_2(s)\|_{W_g'}^2\,ds\\
\leq &\, |\theta_{0m}|_g^2+ \frac{2}{\kappa\Gamma(\alpha)}\int_0^t\|f_2(s)\|_{W_g'}^{2/\alpha_2}\,ds+\frac{2}{\kappa\Gamma(\alpha)}\int_0^t(t-s)^{\frac{\alpha-1}{1-\alpha_2}}\,ds\\
\leq &\, |\theta_{0m}|_g^2+ \frac{2}{\kappa\Gamma(\alpha)}\int_0^T\|f_2(s)\|_{W_g'}^{2/\alpha_2}\,ds+C_2
\end{align*}
where $b_2=\frac{\alpha-1}{1-\alpha_2}$ and $C_2=\frac{2T^{1+b_2}}{\kappa(1+b_2)\Gamma(\alpha)}$. It follows that 
\begin{equation}
\int_0^t (t-s)^{\alpha-1}\|\theta^{(m)}(s)\|_g^2\, ds\leq\, \frac{\Gamma(\alpha)}{\kappa'}|\theta_{0m}|_g^2+\frac{2}{\kappa\kappa'}\int_0^T\|f_2(s)\|_{W_g'}^{2/\alpha_2}\,ds+\frac{\Gamma(\alpha)}{\kappa'}C_2
\end{equation}
On the other hand, integrating  \eqref{3.9} from $0$ to $T$, in the fractional sense, we obtain
\begin{align*}
|u^{(m)}(t)|_g^2 &+\frac{\nu'}{\Gamma(\alpha)}\int_0^t (t-s)^{\alpha-1}\|u^{(m)}(s)\|_g^2\, ds\\
\leq&\, |u_{0m}|_g^2+\frac{c'}{\nu\Gamma(\alpha)}\int_0^t (t-s)^{\alpha-1}\|\theta^{(m)}(s)\|_g^2\, ds+\frac{4}{\nu\Gamma(\alpha)}\int_0^t(t-s)^{\alpha-1}\| f_1(s)\|_{V_g'}^2\,ds\\
\leq &\,|u_{0m}|_g^2+\frac{c'}{\nu\kappa'}|\theta_{0m}|_g^2+ \frac{2c'}{\nu\kappa\kappa'\Gamma(\alpha)}\int_0^t\|f_2(s)\|_{W_g'}^{2/\alpha_2}\,ds+\frac{c'}{\nu\kappa'}C_2\\&+\frac{4}{\nu\Gamma(\alpha)}\int_0^t\|f_1(s)\|_{V_g'}^{2/\alpha_1}\,ds+\frac{4}{\nu\Gamma(\alpha)}\int_0^t(t-s)^{\frac{\alpha-1}{1-\alpha_1}}\,ds\\
\leq &\,|u_{0m}|_g^2+\frac{c'}{\nu\kappa'}|\theta_{0m}|_g^2+ \frac{2c'}{\nu\kappa\kappa'\Gamma(\alpha)}\int_0^t\|f_2(s)\|_{W_g'}^{2/\alpha_2}\,ds+\frac{4}{\nu\Gamma(\alpha)}\int_0^t\|f_1(s)\|_{V_g'}^{2/\alpha_1}\,ds\\
&+C_1
\end{align*}
where $b_1=\frac{\alpha-1}{1-\alpha_1}$ and $C_1=\frac{c'}{\nu\kappa'}C_2+\frac{4T^{1+b_1}}{\nu(1+b_1)\Gamma(\alpha)}$. By using the fact that 
\begin{equation}
\int_0^t (t-s)^{\alpha-1}\|u^{(m)}(s)\|_g^2\, ds\geq T^{\alpha-1}\int_0^t\|u^{(m)}(s)\|_g^2\, ds
\end{equation}and similarly
\begin{equation}
\int_0^t (t-s)^{\alpha-1}\|\theta^{(m)}(s)\|_g^2\, ds\geq T^{\alpha-1}\int_0^t\|\theta^{(m)}(s)\|_g^2\, ds
\end{equation}it follows that
\begin{align}
|u^{(m)}(t)|_g^2 +\frac{\nu' T^{\alpha-1}}{\Gamma(\alpha)}&\int_0^t\|u^{(m)}(s)\|_g^2\, ds\leq
|u_{0m}|_g^2+\frac{c'}{\nu\kappa'}|\theta_{0m}|_g^2\\
&+ \frac{2c'}{\nu\kappa\kappa'\Gamma(\alpha)}\int_0^T\|f_2(s)\|_{W_g'}^{2/\alpha_2}\,ds+\frac{4}{\nu\Gamma(\alpha)}\int_0^T\|f_1(s)\|_{V_g'}^{2/\alpha_1}\,ds +C_1\nonumber
\end{align}

\begin{equation}
|\theta^{(m)}(t)|_g^2 +\frac{\kappa'T^{\alpha-1}}{\Gamma(\alpha)}\int_0^t\|\theta^{(m)}(s)\|_g^2\, ds\leq |\theta_{0m}|_g^2+ \frac{2}{\kappa\Gamma(\alpha)}\int_0^T\|f_2(s)\|_{V_g'}^{2/\alpha_2}\,ds+C_2
\end{equation} 
Consequently
\begin{align}
\displaystyle\sup_{t\in[0,T]}|u^{(m)}(t)|_g^2 \leq
|u_{0m}|_g^2+&\frac{c'}{\nu\kappa'}|\theta_{0m}|_g^2+ \frac{2c'}{\nu\kappa\kappa'\Gamma(\alpha)}\int_0^T\|f_2(s)\|_{W_g'}^{2/\alpha_2}\,ds\\&+\frac{4}{\nu\Gamma(\alpha)}\int_0^T\|f_1(s)\|_{V_g'}^{2/\alpha_1}\,ds +C_1\nonumber
\end{align}

\begin{equation}\label{3.12}
 \displaystyle \sup_{t\in[0,T]}| \theta^{(m)}(t) |^2_g
\leq |\theta_{0m}|_g^2+ \frac{2}{\kappa\Gamma(\alpha)}\int_0^T\|f_2(s)\|_{V_g'}^{2/\alpha_2}\,ds+C_2
\end{equation}
which imply that the sequences $\{u^{(m)}\}_m$ and $\{\theta^{(m)}\}_m$ remain in a bounded set of
$L^{\infty}(0, T ; H_g)$ and $L^{\infty}(0, T ; L^2(\Omega, g))$ respectively. Moreover, for $t=T$, one obtain
\begin{align}\label{u}
\int_0^T\|u^{(m)}(s)\|_g^2\, ds\leq&
\frac{\Gamma(\alpha)}{\nu' T^{\alpha-1}}|u_{0m}|_g^2+\frac{c'}{\nu\nu'\kappa'T^{\alpha-1}}|\theta_{0m}|_g^2+ \frac{2c'}{\nu\nu'\kappa\kappa'T^{\alpha-1}}\int_0^T\|f_2(s)\|_{W_g'}^{2/\alpha_2}\,ds\\
&+\frac{4}{\nu\nu'T^{\alpha-1}}\int_0^T\|f_1(s)\|_{V_g'}^{2/\alpha_1}\,ds +\frac{\Gamma(\alpha)}{\nu' T^{\alpha-1}}C_1\nonumber
\end{align}

\begin{equation}\label{theta}
\int_0^T\|\theta^{(m)}(s)\|_g^2\, ds\leq \frac{\Gamma(\alpha)}{\kappa'T^{\alpha-1}}|\theta_{0m}|_g^2+ \frac{2}{\kappa\kappa'T^{\alpha-1}}\int_0^T\|f_2(s)\|_{V_g'}^{2/\alpha_2}\,ds+\frac{\Gamma(\alpha)}{\kappa'T^{\alpha-1}}C_2
\end{equation} 
which imply that the sequences  $\{u^{(m)}\}_m$ and  $\{\theta^{(m)}\}_m$ remain in a bounded set of
$L^2(0, T ; V_g)$ and $L^2(0, T ; W_g)$ respectively. Consequently, we can assert the existence of elements
$u \in  L^2(0, T ; V_g) \cap L^{\infty}(0, T ; H_g)$ and $\theta \in L^2(0, T ; W_g) \cap L^{\infty}(0, T ; L^2(\Omega, g))$ and the subsequences  $\{u^{(m)}\}_m$ and  $\{\theta^{(m)}\}_m$ such that $u^{(m)} \rightarrow u \in
L^2(0, T ; V_g)$ and $\theta^{(m)} \rightarrow\theta \in L^2(0, T ; W_g)$ weakly and $u^{(m)}\rightarrow u \in L^{\infty}(0, T ; H_g)$ and
$\theta^{(m)} \rightarrow \theta\in  L^{\infty}(0, T ; L^2(\Omega, g))$ weak-star
 as $m \rightarrow\infty$.

Let $\widetilde{u}^{(m)} : \mathbb{R} \rightarrow V_g$
 and $\widetilde{\theta}^{(m)} : \mathbb{R} \rightarrow W_g$ defined as

\[   \widetilde{u}^{(m)}(t)=\left\{
  \begin{array}{rcl}
u^{(m)}(t),&0\leq t\leq T\\
0,& \text{otherwise}
    \end{array} \right.~~~~~\text{and}~~~~~~
      \widetilde{\theta}^{(m)}(t)=\left\{     \begin{array}{rcl}
\theta^{(m)}(t),&0\leq t\leq T\\
0,& \text{otherwise}
    \end{array} \right.     \]
and their Fourier transforms denoted by $\widehat{u}^{(m)}$ and  $\widehat{\theta}^{(m)}$, respectively. We show that the sequence $\{\tilde u^{(m)}\}_m$ remains bounded in $W^\gamma(\mathbb R, V_g, H_g)$ and the sequence $\{\tilde \theta^{(m)}\}_m$ remains bounded in $W^\gamma(\mathbb R, W_g, L^2(\Omega,g))$. To do so, we need to verify that 
\begin{equation}\label{eqgamma1}
\int_{-\infty}^{+\infty}|\tau|^{2\gamma}|\widehat u^{(m)}(\tau)|^2\,d\tau\leq \mathrm{const.}\quad\text{for some }\gamma>0
\end{equation}and
\begin{equation}\label{eqgamma2}
\int_{-\infty}^{+\infty}|\tau|^{2\gamma}|\widehat \theta^{(m)}(\tau)|^2\,d\tau\leq \mathrm{const.}\quad\text{for some }\gamma>0
\end{equation}
In order to prove \eqref{eqgamma1} and \eqref{eqgamma2}, we observe that 
\begin{equation}\label{Fmu}
(D_t^\alpha\tilde u^{(m)},u_k)_g=(\widetilde F^u_m,u_k)_g+(u_{m0},u_k)_g\,I_{-,t}^{1-\alpha}\delta_0-( u^{(m)}(T),u_k)_g \,I_{-,t}^{1-\alpha}\delta_T
\end{equation}
\begin{equation}\label{Fmtheta}
(D_t^\alpha\tilde \theta^{(m)},\theta_k)_g=(\widetilde F^\theta_m,\theta_k)_g+(\theta_{m0},\theta_k)_g\,I_{-,t}^{1-\alpha}\delta_0-( \theta^{(m)}(T),\theta_k)_g \,I_{-,t}^{1-\alpha}\delta_T
\end{equation}
where $\delta_0,\,\delta_T$ are Dirac distributions at $0$ and $T$ and $F_m^u$ and $F_m^\theta$ are defined by
\[
\begin{array}{rcl}
F_m^u=&   \xi \theta^{(m)}+f_1-B_g(u^{(m)},u^{(m)})-\nu A_g u^{(m)}-\nu C_gu^{(m)}\\
F_m^\theta=& f_2-\widetilde{B}_g(u^{(m)},\theta^{(m)}   )-\kappa \widetilde{A}_g\theta^{(m)}+\kappa\widetilde{C}_g\theta^{(m)}+\kappa \widetilde{D}_g\theta^{(m)}
\end{array}
\]for $k = 1,\cdot,\cdot\cdot, m$. Here $\widetilde F_m$ is defined as usual by

\begin{equation}
 \widetilde F_m(t)=
 \left \{
   \begin{array}{r c l}
    F_m(t),&\quad 0\leq t\leq T   \\
    0,&\quad \text{ otherwise}
   \end{array}
   \right .
\end{equation}

Indeed, it is classical that since $\widetilde u^{(m)}$ and $\widetilde \theta^{(m)}$ have two discontinuities at $0$ and $T$, the Caputo derivative of $\widetilde u^{(m)}$ is given by

\begin{align}
D_{-,t}^\alpha\tilde u^{(m)}&= \,I_{-,t}^{1-\alpha}\left(\frac{d}{dt}\tilde u^{(m)}\right)\\
&=\,I_{-,t}^{1-\alpha}\left(\frac{d}{dt} u^{(m)}+u^{(m)}(0)\delta_0-u^{(m)}(T)\delta_T\right)\\
&=\,D_t^\alpha u^{(m)}+I_{-,t}^{1-\alpha}\left(u^{(m)}(0)\delta_0-u^{(m)}(T)\delta_T\right)
\end{align}and the one of $\widetilde \theta^{(m)}$ is given by
\begin{equation}
D_{-,t}^\alpha\widetilde \theta^{(m)}=\,D_t^\alpha \theta^{(m)}+I_{-,t}^{1-\alpha}\left(\theta^{(m)}(0)\delta_0-\theta^{(m)}(T)\delta_T\right)
\end{equation}

By the Fourier transform, \eqref{Fmu} and \eqref{Fmtheta} yield to
\begin{align}\label{fourieru}
(2i\pi\tau)^\alpha(\widehat u^{(m)},u_k)_g=&(\widehat{F}^u_m,u_k)_g+(u_{m0},u_k)_g(2i\pi\tau)^{\alpha-1}\\
&\qquad\qquad-(u^{(m)}(T),u_k)_g(2i\pi\tau)^{\alpha-1}e^{-2i\pi T\tau}
\end{align}
\begin{align}\label{fouriertheta}
(2i\pi\tau)^\alpha(\widehat \theta^{(m)},\theta_k)_g=&(\widehat{F}^\theta_m,\theta_k)_g+(\theta_{m0},\theta_k)_g(2i\pi\tau)^{\alpha-1}\\
&\qquad\qquad-(\theta^{(m)}(T),\theta_k)_g(2i\pi\tau)^{\alpha-1}e^{-2i\pi T\tau}
\end{align}here $\widehat u^{(m)}$ and $\widehat{F}_m$ denote the Fourier transforms of $\tilde u^{(m)}$ and $\widetilde{F}_m$, respectively. We multiply \eqref{fourieru} and \eqref{fouriertheta} by $\widehat{f}_j^{(m)}$ and $\widehat{g}_j^{(m)}$ respectively and add
these equations for $k= 1,\dots, m$ to get
\begin{align}
(2i\pi\tau)^\alpha|\widehat u^{(m)}(\tau)|_g^2=&(\widehat{F}^u_m(\tau),\widehat u^{(m)}(\tau))_g+(u_{m0},\widehat u^{(m)}(\tau))_g(2i\pi\tau)^{\alpha-1}\\
&\qquad\qquad-(u^{(m)}(T),\widehat u^{(m)}(\tau))_g(2i\pi\tau)^{\alpha-1}e^{-2i\pi T\tau}
\end{align}
\begin{align}
(2i\pi\tau)^\alpha|\widehat \theta^{(m)}(\tau)|_g^2=&(\widehat{F}^\theta_m(\tau),\widehat \theta^{(m)}(\tau))_g+(\theta_{m0},\widehat \theta^{(m)}(\tau))_g(2i\pi\tau)^{\alpha-1}\\
&\qquad\qquad-(\theta^{(m)}(T),\widehat \theta^{(m)}(\tau))_g(2i\pi\tau)^{\alpha-1}e^{-2i\pi T\tau}
\end{align}
Since the integrals on the right hand side of the inequalities
\begin{align}
\displaystyle \int^T_0 \|F_m^u(t) \|_{V'_g}dt\leq \int^T_0 c( |\xi|_{\infty} \| \theta^{(m)}(t) \|_g&+\|f_1(t)\|_{V'_g}+| u^{(m)}(t)|_g \|u^{(m)} \|_g\\
&+ \|u^{(m)}(t)  \|_g+| \nabla g |_{\infty}\|u^{(m)}(t)  \|_g)dt\nonumber
\end{align}
\begin{align}
\displaystyle \int^T_0 \|F_m^\theta(t) \|_{W'_g}dt\leq \int^T_0 c'(  \| f_2(t) \|_{W'_g}+& |u^{(m)}(T)|_g \|\theta^{(m)} (t) \|_g+ \|\theta^{(m)} (t) \|_g\\
&  +|\nabla g|_{\infty} \| \theta^{(m)}(t) \|_g+|\Delta g|_{\infty}\| \theta^{(m)}(t) \|_g)dt\nonumber
\end{align}remains bounded, $ \|F_1(t) \|_{V'_g}$ and $\|F_2(t) \|_{W'_g}$ are bounded in $L^1(0, T ; V'_g)$ and $L^1(0, T ; W'_g)$ respectively. Therefore, for all $m$

\[
\begin{array}{ ccc}
\displaystyle \sup_{\tau\in\mathbb{R}} \| \widehat{F}^u_m(\tau) \|_{V'_g} \leq c_1&\text{and}&  \displaystyle \sup_{\tau\in\mathbb{R}} \| \widehat{F}^\theta_m(\tau) \|_{W'_g}\leq c_2.
\end{array}\]
Moreover, since $u^{(m)}(0)$,  $u^{(m)}(T),$  $\theta^{(m)}(0)$ and $\theta^{(m)}(T )$ are bounded, we get 

\begin{align*}
    |\tau|^\alpha | \widetilde{u}^{(m)}(\tau)|^2_g      \leq &\, c_1\| u^{(m)} \|_{V_g}+c_2 |\tau|^{\alpha-1}\,|u^{(m)}  |_g\\
    \leq&\, c_3 \| u^{(m)} \|_{V_g}
    \end{align*}
\begin{align}
    |\tau|^\alpha | \widetilde{\theta}^{(m)}(\tau)|^2_g      \leq &\, c'_1\| \theta^{(m)} \|_{W_g}+c'_2 |\tau|^{\alpha-1}\,|\theta^{(m)}  |_g\\
    \leq &\,c_3 \| \theta^{(m)} \|_{W_g}
\end{align}
For $\gamma$ fixed, $\gamma<\alpha/4$, we observe that 
\[
|\tau|^{2\gamma}\leq c(\gamma)\frac{1+|\tau|^\alpha}{1+|\tau|^{\alpha-2\gamma}}
\]
Then we can write
\begin{align*}
\int_{-\infty}^{+\infty}|\tau|^{2\gamma}|\widehat u^{(m)}(\tau)|^2_g\leq & c_5(\gamma) \int_{-\infty}^{+\infty}\frac{1+|\tau|^\alpha}{1+|\tau|^{\alpha-2\gamma}}|\widehat u^{(m)}(\tau)|^2_g\,d\tau\\
\leq & c_6(\gamma) \int_{-\infty}^{+\infty}\frac{1}{1+|\tau|^{\alpha-2\gamma}}\|\widehat u^{(m)}(\tau)\|^2_{V_g}\,d\tau\\
&\qquad\qquad+c_7(\gamma) \int_{-\infty}^{+\infty}\frac{|\tau|^{\alpha-1}}{1+|\tau|^{\alpha-2\gamma}}\|\widehat u^{(m)}(\tau)\|^2_{V_g}\,d\tau
\end{align*}
\begin{align*}
\int_{-\infty}^{+\infty}|\tau|^{2\gamma}|\widehat \theta^{(m)}(\tau)|^2_g
\leq & c'_6(\gamma) \int_{-\infty}^{+\infty}\frac{1}{1+|\tau|^{\alpha-2\gamma}}\|\widehat \theta^{(m)}(\tau)\|^2_{W_g}\,d\tau\\
&\qquad\qquad+c'_7(\gamma) \int_{-\infty}^{+\infty}\frac{|\tau|^{\alpha-1}}{1+|\tau|^{\alpha-2\gamma}}\|\widehat \theta^{(m)}(\tau)\|^2_{W_g}\,d\tau
\end{align*}
By Parseval inequality, the first integral is bounded as $m\to\infty$. Applying the Schwartz inequality, the second integrals yield to
\begin{align}
\int_{-\infty}^{+\infty}\frac{|\tau|^{\alpha-1}}{1+|\tau|^{\alpha-2\gamma}}\|\hat u^{(m)}(\tau)\|^2_g\,d\tau\leq&\left(\int_{-\infty}^{+\infty}\frac{d\tau}{(1+|\tau|^{\alpha-2\gamma})^2}\right)^{1/2}\\
&\qquad\qquad\times\left(\int_{-\infty}^{+\infty}|\tau|^{2\alpha-2}\|\hat u^{(m)}(\tau)\|^2_g\,d\tau\right)^{1/2}
\end{align}
\begin{align}
\int_{-\infty}^{+\infty}\frac{|\tau|^{\alpha-1}}{1+|\tau|^{\alpha-2\gamma}}\|\widehat \theta^{(m)}(\tau)\|^2_g\,d\tau\leq&\left(\int_{-\infty}^{+\infty}\frac{d\tau}{(1+|\tau|^{\alpha-2\gamma})^2}\right)^{1/2}\\
&\qquad\qquad\times\left(\int_{-\infty}^{+\infty}|\tau|^{2\alpha-2}\|\widehat \theta^{(m)}(\tau)\|^2_g\,d\tau\right)^{1/2}
\end{align}
The first integrals are finite due to $\gamma<\alpha/4$. On the other hand, it follows from the Parseval equality that
\begin{align*}
\int_{-\infty}^{+\infty}|\tau|^{2\alpha-2}\|\hat u^{(m)}(\tau)\|^2_g\,d\tau&=\int_{-\infty}^{+\infty}\|\,_{-\infty}\mathrm I_t^{1-\alpha}\tilde u^{(m)}(t)\|_g^2\,dt\\
&=\int_0^T\|\,_0\mathrm I_t^{1-\alpha} u^{(m)}(t)\|^2_g\,dt\\
&\leq\left(\frac{T^{1-\alpha}}{\Gamma(2-\alpha)}\right)^2\int_0^T\|u^{(m)}(t)\|_{V_g}^2\,dt
\end{align*}
\begin{align*}
\int_{-\infty}^{+\infty}|\tau|^{2\alpha-2}\|\widehat \theta^{(m)}(\tau)\|^2_g\,d\tau\leq\left(\frac{T^{1-\alpha}}{\Gamma(2-\alpha)}\right)^2\int_0^T\|\theta^{(m)}(t)\|_{W_g}^2\,dt
\end{align*}

Which implies that \eqref{eqgamma1} and \eqref{eqgamma2} holds. We know that subsequence of $\{u^{(m)}\}_m$ and $\{\theta^{(m)}\}_m$ (which we will denote with the same symbol) converge to some $u$ weakly in $L^2(0,T;V_g)$ and weak-star in $L^\infty(0,T;H_g)$ with $u\in L^2(0,T;V_g)\cap L^\infty(0,T;H_g)$. Similarly there exists a subsequence of $\{\theta^{(m)}\}_m$ (which we will denote with the same symbol) converge to some $\theta$ weakly in $L^2(0,T;W_g)$ and weak-star in $L^\infty(0,T;L^2(\Omega,g))$ with $\theta\in L^2(0,T;W_g)\cap L^\infty(0,T;L^2(\Omega,g))$. As $W^\gamma(0,T,V_g;H_g)$ is compactly embedded in $L^2(0,T; H_g)$ and $W^{\gamma}(\mathbb{R}, W_g, L^2(\Omega,g))$ in $L^2(0, T ; L^2(\Omega, g))$ then $\{u^{(m)}\}_m$ strongly converges in $L^2(0,T;H_g)$ and  $\{\theta^{(m)}\}_m$ in $L^2(0, T ; L^2(\Omega, g))$ respectively.

In order to pass to the limit, we consider the scalar functions $\Psi_1(t)$ and $\Psi_2(t)$ continuously differentiable on $[0, T ]$ and such that 
$\Psi_1(T)=0$ ve $\Psi_2(T)=0$. We multiply (\ref{3.4})
and (\ref{3.5}) by $\Psi_1(t)$ and $\Psi_2(t)$ respectively and then integrate by parts. This leads to the equations
\begin{align*}
\int_0^T&(u^{(m)}(t),\mathrm D_{t,T}^\alpha\Psi_1(t)u_k)_g\,dt+\int_0^T b_g(u^{(m)}(t),u^{(m)}(t),\Psi_1 u_k)\,dt\\
+&\nu\int_0^T((u^{(m)}(t),\Psi_1 u_k))_g+\nu\int_0^Tb_g(\frac{\nabla g}{g},u^{(m)}(t),\Psi_1 u_k)\,dt=(u_{0m},\mathrm I_{0,T}^{1-\alpha}\Psi_2(t) u_k)_g\\
+&\int^T_0 (\xi\theta^{(m)}(t),\Psi_1u_k)_gdt+\int_0^T(f_1(t),u_k)_g\,dt
\end{align*}

\begin{align*}
\int_0^T&(\theta^{(m)}(t),\mathrm D_{t,T}^\alpha\Psi_2(t)\theta_k)_g\,dt+\int_0^T \widetilde b_g(u^{(m)}(t),\theta^{(m)}(t),\Psi_2 \theta_k)\,dt\\
+&\kappa\int_0^T((\theta^{(m)}(t),\Psi_2 \theta_k))_g\,dt+\kappa\int_0^T\widetilde b_g(\frac{\nabla g}{g},\theta_k,\Psi_2 \theta^{(m)}(t))\,dt=(\theta_{0m},\mathrm I_{0,T}^{1-\alpha}\Psi_2(t) \theta_k)_g\\
+&\int_0^T(f_2(t),\Psi_2\theta_k)_g\,dt
\end{align*}
Following the technique given in \cite{G20,T84}, as $m \rightarrow \infty$ we obtain
\begin{align}\label{eqproof1}
\int_0^T&(u(t),\mathrm D_{t,T}^\alpha\Psi_1(t)u_k)_g\,dt+\int_0^T b_g(u(t),u(t),\Psi_1 u_k)\,dt+\nu\int_0^T((u(t),\Psi_1 u_k))_g\\
&+\nu\int_0^Tb_g(\frac{\nabla g}{g},u(t),\Psi_1 u_k)\,dt=(u_{0},\,\mathrm I_{0,T}^{1-\alpha}\Psi_1 u_k)_g+\int^T_0 (\xi\theta(t),\Psi_1v)_gdt+\int_0^T(f_1(t),u_k)_g\,dt
\end{align}

\begin{align}\label{eqproof2}
\int_0^T&(\theta(t),\mathrm D_{t,T}^\alpha\Psi_2(t)\theta_k)_g\,dt+\int_0^T \widetilde b_g(u(t),\theta(t),\Psi_2 \theta_k)\,dt\\\nonumber
+&\kappa\int_0^T((\theta(t),\Psi_2 \theta_k))_g\,dt+\kappa\int_0^T\widetilde b_g(\frac{\nabla g}{g},\theta_k,\Psi_2 \theta(t))\,dt=(\theta_{0},\mathrm I_{0,T}^{1-\alpha}\Psi_2(t) \theta_k)_g\\\nonumber
+&\int_0^T(f_2(t),\Psi_2\theta_k)_g\,dt
\end{align}
This equations hold for $v$ and $\tau$ which are finite linear combination of $u_k$ and $\theta_k$, $k=1,\dots,m$ and by continuity it holds for any $v$ in $V_g$ and $\tau\in H_g$. It then follows that $\{u,\theta\}$ satisfies the two first equations \eqref{benard}. To end the proof it still to check that $\{u,\theta\}$ satisfies the initial conditions $u(0)=u_0$ and $\theta(0)=\theta_0$. To do so it suffices to multiply the two first equations in \eqref{benard} by $\Psi_1$ and \eqref{} by $\Psi_2$ respectively and then integrate. By making use of the integration by parts and comparing with \eqref{eqproof1} and \eqref{eqproof2}, one can find that 
\[
(u_{0}-u(0), v)_g\,\mathrm I_{0,T}^{1-\alpha}\Psi_2(t)=0,\quad\text{and}\quad (\theta_{0}-\theta(0), \tau)_g\,\mathrm I_{0,T}^{1-\alpha}\Psi_2(t)=0
\]which lead to the desired result by taking a particular choice of $\Psi_1$ and $\psi_2$.

%%%%%%%%%%%%%%%%

For the uniqueness of the weak solutions let $(u_1,\theta_1)$ and $(u_2,\theta_2)$ be two weak solutions
with the same initial condition. Let $w = u_1- u_2$ and $\widetilde{w}~ = \theta_1 -\theta_2$. Then we have

\[
\begin{array}{r}
D^{\alpha}_{t}(w,v)_g+ b_g(u_1,u_1,v)- b_g(u_2,u_2,v)+\nu (\nabla w,\nabla v  )_g+\nu (C_gw,v)_g= (\xi  \widetilde{w},v)_g
\\  \\ D^{\alpha}_{t} ( \widetilde{w},\tau)_g+ \widetilde{b}_g(u_1,\theta_1,\tau)-\widetilde{b}_g(u_2,\theta_2,\tau)    +\kappa(  \nabla\widetilde{w},\nabla\tau )_g+\kappa \widetilde{b}_g
( \displaystyle \frac{\nabla g}{g},\tau,\widetilde{w})=0
\end{array}
\]
Taking $v = w(t)$ and $\tau = \widetilde{ w}(t)$ one obtains,
\[  \displaystyle D^{\alpha}_{t}(w,w)_g+b_g(w,u_2,w)+\nu | A^{1/2}_g w |^2_g+\nu (C_g w,w)_g=(\xi\widetilde{w},w)g  \]
\[  \displaystyle D^{\alpha}_{t} ( \widetilde{w},\widetilde w)_g+\widetilde{b}_g(u_1,\theta_1,\widetilde{w})-\widetilde{b}_g(u_2,\theta_2,\widetilde{w})
+
\kappa| \widetilde{A}^{1/2}_g \widetilde{w} |^2_g+ \kappa \widetilde{b}_g ( \displaystyle \frac{\nabla g}{g},\widetilde{w},\widetilde{w} )=0.\]
By applying the bounds on the terms $b_g$,  $\widetilde{b}_g$ it then follows by Cauchy-Schwarz inequality and Gronwall like inequality that $w(t) = 0$ and $\widetilde{w}(t) = 0$ for all $t \geq 0$ since we have $w(0) = 0$
and $\widetilde{w}~(0) = 0$. Thus the Theorem is proved.

\end{proof}

\end{document}